\documentclass[12pt]{article}
\usepackage{amsmath}
\usepackage{amsthm}
\usepackage{amsfonts}
\usepackage{amssymb}
\usepackage{latexsym} 
\usepackage[toc,page]{appendix}
\usepackage{dsfont}
\usepackage{microtype}

\usepackage[matrix,tips,graph,curve]{xy}



\makeatletter 
\@addtoreset{equation}{section}
\makeatother

\theoremstyle{plain}
\newtheorem{theorem}[equation]{Theorem}
\newtheorem{corollary}[equation]{Corollary}

\newtheorem{proposition}[equation]{Proposition}

\theoremstyle{definition}
\newtheorem{definition}[equation]{Definition}

\newtheorem{remark}[equation]{Remark}

\newcommand{\IC}{\mathbb{C}}

\newcommand{\IQ}{\mathbb{Q}}
\newcommand{\IR}{\mathbb{R}}
\newcommand{\IS}{\mathbb{S}}

\newcommand{\IZ}{\mathbb{Z}}

\newcommand{\End}{\mathrm{End}}

\newcommand{\Hom}{\mathrm{Hom}}
\newcommand{\Aut}{\mathrm{Aut}}

\newcommand{\rank}{\mathrm{rank}}
\renewcommand\dim{{\rm dim\,}}

\def\d/{/\mspace{-6.0mu}/}

\newcommand{\mf}{\mathfrak}

\usepackage[
style=alphabetic]{biblatex}
\usepackage[utf8]{inputenc}
 \usepackage[T1]{fontenc} 
\usepackage{longtable}
\usepackage{breakurl}
\usepackage[breaklinks]{hyperref}

\usepackage{commath}

\begin{document}
\title{Hodge Representations of Calabi-Yau 3-fold Type}
\author{Xiayimei Han}
\date{}
\maketitle
{\footnotesize{\noindent Dedicated to family, friends and mentors for generous support, and Andre for good company.}}
\vspace{1cm}

\begin{abstract}
\noindent In this thesis, I apply the Green-Griffiths-Kerr classification of Hodge representations to enumerate the Lie algebra Hodge representations of CY 3-fold type.
\end{abstract}

\vspace{0.5cm}
\section{Introduction}
In this thesis, we will complete the classification of real Hodge representations of Calabi-Yau 3-fold type that was begun in \cite{Friedman_2013}. 
In the following subsections, we will introduce Hodge representation of Lie groups and Lie algebras, and the Hodge representation of Calabi-Yau 3-fold type. The introductory materials are mostly taken from \cite{10.2307/j.ctt7s3h1}, \cite{Robles2014} and \cite{han_robles_2020}.
\subsection{Hodge Representation on Lie group and Lie algebra Level}
\begin{definition}
An \textbf{algebraic group} is a group that is an algebraic variety, such that the multiplication and inversion operations are given by maps on the variety that are locally given by polynomials.
\end{definition}

\begin{definition}
A {\bf{Hodge structure}} of weight $n\in \IZ$ on a rational vector space $V$ is a non-constant homomorphism $\phi:S^1\to \Aut(V_\IR)$ of $\IR$-algebraic groups such that $\phi(-1)=(-1)^n\mathds{1}$. The associated {\bf Hodge decomposition} $V_\IC=\oplus_{p+q=n} V^{p,q}$ with $\overline{V^{p,q}}=V^{q,p}$ is given by $$V^{p,q}=\{v\in V_\IC| \phi(z)v=z^{p-q}v\quad \forall z\in S^1\}.$$ The {\bf Hodge numbers} $h=(h^{p,q})$ are $h^{p,q}=\dim_\IC V^{p,q}$.
\end{definition}

\begin{definition}
Let $Q:V_\IC \times V_\IC \to \IQ$ be a non-degenerate bilinear form satisfying $Q(u,v)=(-1)^nQ(v,u)$ for all $u,v\in V_\IC$. A Hodge structure $(V,\phi)$ of weight $n$ is {\textbf{polarized}} by $Q$ if the {\textbf{Hodge-Riemann bilinear relations}} hold:
$$Q(V^{p,q},V^{r,s})=0\quad\text{unless }q=r,p=s,$$
$$Q(v,\phi(i)\bar{v})>0\quad \forall 0\neq v\in V_\IC.$$
Throughout we assume the Hodge structure $\phi$ is polarized.
\end{definition}
Finally, to introduce the motivation of this thesis, we need to define Mumford-Tate groups and Mumford-Tate domains. To do this, we follow the construction in \cite{MTD}.

\begin{definition}
Let $V$ be a rational vector space and $V\otimes \IC=V_\IC=\oplus_{p+q=n} V^{p,q}$ with $\overline{V^{p,q}}=V^{q,p}$ be a Hodge decomposition. Use $\IS(\IR)$ to denote the matrix group $\{\begin{pmatrix}
a & b\\
-b & a
\end{pmatrix}\in \text{GL}_2\IR|a,b\in \IR\},$ and use $M(a,b)$ to denote $\begin{pmatrix}
a & b\\
-b & a
\end{pmatrix}\in \IS(\IR)$. Then we can define an algebraic group representation $\psi:\IS(\IR)\to \Aut(V_\IC)$ by letting $M(a,b)$ act on $H^{p,q}$ as multiplication by $(a+ib)^p(a-ib)^q$. Then the \textbf{Mumford-Tate group $M_\psi$} of the Hodge decomposition is the smallest algebraic subgroup of $GL(V)$ defined over $\IQ$, whose set of real points contains $\psi(S^1)$, where we identify $S^1$ with $\{\begin{pmatrix}
a & b\\
-b & a
\end{pmatrix}:a^2+b^2=1\}\subset \IS(\IR)$. The \textbf{Mumford-Tate domain} is the $M_\psi(\IR)$ orbit of $\psi$.
\end{definition}

In order to classify Mumford-Tate groups and Mumford-Tate domains, Griffiths, Green, and Kerr introduced the following notion:
\begin{definition}
A {\bf Lie group Hodge representation $
(G,\rho,\phi)$} consists of a reductive $\IQ$-algebraic group $G$, a faithful representation $\rho:G\to \Aut(V,Q)$, defined over $\IQ$ and a non-constant homomorphism $\phi:S^1\to G_\IR$ of $\IR$-algebraic groups such that $(V,Q,\rho\circ\phi)$ is a polarized Hodge structure. Note that $G_\IR$ denotes the set of real points of $G$.
\end{definition}

Some results on the classification of Lie group Hodge representations can be seen in \cite{10.2307/j.ctt7s3h1} and \cite{Zarhin1983}. To avoid the complexities of working with $\IQ$-algebraic groups, we "forget" the rational structure and work with the following notion:

\begin{definition}
A {\bf real Lie group Hodge representation $(G,\rho, \phi)$} is as in Definition 1.4, but $G$ a reductive $\IR$-algebraic group and $V$ a real vector space. Such $G$ is called a {\bf real Hodge group}.
\end{definition}

Now we proceed to define real Lie algebra Hodge representations. To do this, we first cite the following fundamental result from Section IV of \cite{10.2307/j.ctt7s3h1}:

\begin{theorem}\label{GGK thm}
Assume $G$ is a semisimple $\IQ$-algebraic group. If $G$ has a Hodge representation, then $G_\IR$ contains a compact maximal torus $T$ with $\phi(S^1)\subset T$ and $\dim(T)=\rank(G)$. Moreover, the maximal subgroup $H_\phi\subset G_\IR$ such that $\rho(H_\phi)$ fixes the polarized Hodge structure $(V,Q,\rho\circ \phi)$ is compact and is the centralizer of the circle $\phi(S^1)$ in $G_\IR$. 
\end{theorem}

Let $G$ be a Hodge group. Then $G$ is a Lie group and we denote its Lie algebra by $\mf{g}$. Fixing a Cartan subalgebra $\mf{t}$, we know that $\mf{g}$ has corresponding weight lattice $\Lambda_w$ and root lattice $\Lambda_R$. A grading element is an element $E$ of $\mf{t}$ such that $\lambda(E)\in \IZ$ for all $\lambda\in \Lambda_R$. Let $\tilde{G}$ be the simply connected simple Lie group with Lie algebra $\mf{g}$. Then its center $Z(\tilde{G})=\frac{\Lambda_w}{\Lambda_R}$. We also know that all connected simple Lie groups with Lie algebra $\mf{g}$ are of the form $\tilde{G}/C$ with $C$ a subgroup of $Z(\tilde{G})$, so the simple Lie groups with Lie algebra $\mf{g}$ are in one to one correspondence with subgroups of $Z(\tilde{G})$, each of which is isomorphic to some $\frac{\Lambda}{\Lambda_R}$ where $\Lambda$ is a sublattice of $\Lambda_w$ containing $\Lambda_R$. Suppose that under this correspondence, $G$ is labeled by $\Lambda$. Let $\Lambda^* :=\Hom (\Lambda,2 \pi i \IZ).$
Then the maximal torus of $G$ is compact by Theorem $\ref{GGK thm}$, and is given by 
\begin{equation*}
    T\simeq \mf{t}/\Lambda^*. 
\end{equation*}

\indent Suppose that $\phi:S^1\to G$ is a homomorphism of real algebraic groups with image $\phi(S^1)\subset T.$ Let $L_\phi$ be the lattice point such that $\phi(e^{2\pi i t})=tL_\phi \mod \Lambda^*$ for $t\in \IR$. Fix a representation $\rho: G \to \Aut(V,Q)$ defined over $\IQ$. Suppose that $\rho \circ \phi$ defines a polarized Hodge structure on $V$. Then the group element $\rho\circ \phi(e^{2\pi i t})\in \Aut(V_\IC)$ acts on $V^{p,q}$ by $e^{2\pi i(p-q)t}\mathds{1}$. Taking the derivative with respect to $t$, we find that $L_\phi$ acts on $V^{p,q}$ by $2\pi i (p-q)\mathds{1}$. Define $T_\phi\in \Hom(\Lambda, \frac{1}{2}\IZ)\subset i\mf{t}\subset \mf{t}_\IC$ by
\begin{equation}\label{grading element}
    T_\phi:=\frac{L_\phi}{4\pi i}.
\end{equation} Then $T_\phi$ acts on $V^{p,q}$ by the scalar $(p-q)/2=p-n/2$. Therefore, the $T_\phi$-graded decomposition of $V_\IC$ is 
$$V_\IC=V_{m/2}\oplus V_{m/2-1}\oplus\cdots \oplus V_{1-m/2}\oplus V_{-m/2},$$
where $0\leq m\leq n$ is the largest integer such that $V_{m/2}\neq 0$ and $V_l=\{v\in V_\IC|T_\phi(v)=lv\}$. In this notation, $V^{p,q}=V_{(p-q)/2}=V_{p-n/2}$ and $F^p_\phi=V_{p-n/2}\oplus V_{p+1-n/2}\oplus\cdots \oplus V_{m/2}$. By this construction, we get the Lie algebra Hodge representation:
\begin{definition}
A {\bf real Lie algebra Hodge representation} $(\mf{g}_\IR, E, V_\IR)$ corresponding to a real Lie group representation $(G,\rho,\phi)$ consists of the Lie algebra $\mf{g}_\IR$ of $G$, grading element $E$ equal to $T_\phi$ in $\eqref{grading element}$, and $V_\IR$ is the underlying real vector space of $V_\IC$. 
\end{definition}
In practice, real Lie algebra Hodge representations are easier to compute than real Lie group Hodge representations and yield important insights into the latter. In this paper, we will classify all the real Hodge Representations of $CY$ 3-fold type.

\subsection{Hodge representations associated to CY 3-fold type}
The goal of this thesis is to identify the real Lie algebra Hodge representations that can arise when $V_\IR$ is the period domain parameterizing polarized Hodge structure of CY 3-fold type. These consist of the data $E\in i\mf{g_\IR}\subset \End(V_\IC,Q)$, where
\begin{enumerate}
    \item $V_\IR$ is a real space. 
    \item $Q:V_\IR \times V_\IR\to \IR$ is a non-degenerate skew form. Thus $\Aut(V_\IR,Q)\simeq Sp_{2d}\IR$, where $2d=\dim V_\IR$.
    \item $\mf{g}_\IR$ is a real reductive Lie subalgebra of $\End(V_\IR,Q)\simeq \mf{sp}_{2d}\IR$.
    \item $E$ is a semisimple element of $i\mf{g}_\IR \subset \mf{g}_\IC$ acting on $V_\IC=V_\IR\otimes \IC$ with eigenspace decomposition $V_\IC=V_{\frac 3 2}\oplus V_{\frac 1 2}\oplus V_{-\frac 1 2}\oplus V_{-\frac 3 2}$, and $\dim V_{\frac{3}{2}}$.  
    \item The complex dimensions of $V_{\frac 3 2},V_{\frac 1 2}, V_{-\frac 1 2},V_{-\frac 3 2}$ respectively are 1, $a$, $a$, 1 for some positive integer $a$ and $\overline{V_l}=V_{-l}$.
\end{enumerate}
In particular, the constraint conditions that the representations are level 3 (see Subsection 1.3 for specific definition of level) and have first Hodge number $h^{3,0}=1$ ensure that we get Hodge representations associated to CY 3-fold type. For more background information, the reader might consult Section 1 and 2 of \cite{han_robles_2020} as well as Appendix B of \cite{Robles2014}. In order to phrase our results, we assume a fixed torus and a fixed Weyl Chamber. To enumerate all such tuples of $(\mf{g}_\IR, E, V_\IR)$, we take the following steps:

\subsection{Step 1: Reduce to the case that \texorpdfstring{$V_\IR$} is irreducible}
Suppose $V_\IR=V^1_\IR\oplus V^2_\IR$ as $\mf{g}_\IR$ representations. Then $V_\IC=V^1_\IC\oplus V^2_\IC$. Hence, $$V_\IC=V_{\frac 3 2}\oplus V_{\frac 1 2}\oplus V_{-\frac 1 2}\oplus V_{-\frac 3 2}$$ if and only if 
$$V^1_\IC=V^1_{\frac 3 2}\oplus V^1_{\frac 1 2}\oplus V^1_{-\frac 1 2}\oplus V^1_{-\frac 3 2}\qquad \text{and }\qquad V^2_\IC=V^2_{\frac 3 2}\oplus V^2_{\frac 1 2}\oplus V^2_{-\frac 1 2}\oplus V^2_{-\frac 3 2}.$$
Moreover, $\dim V_l=\dim V^1_l+\dim V^2_l$ so condition 5 above translates to 
$$1=\dim V_{\pm\frac 3 2}=\dim V^1_{\pm\frac 3 2}+\dim V^2_{\pm\frac 3 2 }.$$
Thus, without loss of generality, $V^1_\IC=V^1_{\frac 1 2}\oplus V^1_{-\frac 1 2}$ and $V^2_\IC=V^2_{\frac 3 2}\oplus V^2_{\frac 1 2}\oplus V^2_{-\frac 1 2} \oplus V^2_{-\frac 3 2}$, where $\dim V^2_{\frac 3 2}=\dim V^2_{-\frac{3}{2}}=1$.  And we say that $V^1_\IR$ is a level 1 Hodge representation and $V^2_\IR$ is a level 3 Hodge representation. These two types of representations are classified respectively in Sections 2 and 3. 

\subsection{Step 2: Reduction to complexification}
To simplify the computations, instead of directly computing $\mf{g}_\IR$ and $V_\IR$, we will first work with the complexifications $\mf{g}_\IC=\mf{g}_\IR\otimes \IC$ and $V_\IC=V_\IR\otimes \IC$, which will be carried out in the following subsections. Then we can give explicit formulae of $\mf{g}_\IR$ from $\mf{g}_\IC$ in most cases, namely when $E=A^i$, the unique Cartan element that is dual to simple root $\alpha_i$ for some $i$, although recovering $V_\IR$ is more complicated. In order to recover $\mf{g}_\IR$, one should first observe that the grading element $E$ corresponds to a noncompact root as explained in Chapter IV of \cite{10.2307/j.ctt7s3h1}. Then by consulting \cite{Knapp}, we enumerate all the underlying real forms $\mf{g}_\IR$ in Sections 2 and 3. 

\subsection{Step 3: Reduction to semisimple}
 Given an irreducible faithful representation $U$ of $\mf{g}_\IC$ with highest weight $\mu$, we have $\mf{g}_\IC\subset \Aut(U)$. Under the adjoint representation, $\mf{g}_\IC$ decomposes as $\mf{g}^{ss}_\IC\oplus \mf{z}$, where $\mf{g}^{ss}_\IC$ denotes the semisimple part and $\mf{z}$ denotes the center. By Schur's lemma (see Appendix C), if the semisimple part of $\mf{g}_\IC$ is simple, then $\mf{z}$ is 0 or 1 dimensional. In general, the dimension of $\mf{z}$ does not exceed the number of simple direct summands in $\mf{g}^{ss}_\IC$.  Given a semisimple element $E$, $E=E_{ss}+E'$ for unique choices of $E_{ss}\in \mf{g}_\IC^{ss}$ and $E'\in\mf{z}$. By virtue of infinitesimal period relation(IPR) on the dual of the Mumford-Tate domain, we without loss of generality assume that $\alpha_i(E_{ss})\in\{0,1\}$ for each simple root $\alpha_i$ throughout this paper \cite{Robles2014}. On $U$, $E_{ss}$ acts via weights and $E'$ acts as $c\mathds{1}$ for some $c\in \IQ$. Thus, $\mu(E)=\mu(E_{ss})+ c$ on $U$ and $\mu^*(E)=\mu^*(E_{ss})- c$ on $U^*$, where $\mu^*$ is the highest weight of $\mf{g}_\IC$ on $U^*$. In order to recover $\mf{g}_\IR$, we need the tuples $(\mf{g}_\IC, E_{ss},c )$, which we enumerate in the following Sections 2 and 3. More specifically, the grading element $E$ determines the maximal compact subalgebra $\mf{k}$ of $\mf{g}_\IR^{ss}$, which is just the direct sum of the even $E$-eigenspaces in the adjoint representation. This is enough to determine the real form $\mf{g}_\IR^{ss}$ using Vogan diagram classification, and \cite{Knapp} is a good source of information. Then if $c\neq 0$, we have $\mf{g}_\IC=\mf{g}_\IC^{ss}\oplus \IC$ and therefore $\mf{g}_\IR=\mf{g}_\IR^{ss}\oplus \IR.$

\subsection{Step 4: Real, Complex, or Quaternionic representations}
When $V_\IR$ is an irreducible representation of $\mf{g}_\IR$, there are three possible cases:
\begin{equation*}
    V_\IR\otimes \IC =\left\{ 
    \begin{array}{ll}
        U\text{ and }U=U^*\text{ as representations of $\mf{g}_\IC$} \\ \hfill \quad (U\text{ is a  \textbf{real} representation of $\mf{g}_\IR$})\\
        U\oplus U^* \text{ and }U\neq U^*\text{ as representations of $\mf{g}_\IC$} \\
        \hfill \quad (U\text{ is a \textbf{complex} representation of $\mf{g}_\IR$})
        \\
        U\oplus U^* \text{ and }U=U^* \text{ as representations of $\mf{g}_\IC$} \\\hfill \quad  (U\text{ is a \textbf{quaternionic} representation of $\mf{g}_\IR$})
    \end{array}
    \right.
\end{equation*}
where $U$ is an irreducible representation of $\mf{g}_\IC$ with highest weight $\mu$.

 To distinguish the three possible cases of $U$, one could immediately observe that once we have $U\neq U^*$ (or equivalently $\mu\neq \mu^*$), $U$ is complex. When $U=U^*$ (or equivalently $\mu=\mu^*$), we determine whether $U$ is real or quaternionic using the following test:
Define $$H_\phi:=2\sum_{\alpha_j(E_{ss})= 0}A^j. $$
Then $U$ is quaternionic if and only if $\mu(H_\phi)$ is odd, and real if and only if $\mu(H_\phi)$ is even. 
\newline In general, when $U$ is real or quaternionic, we necessarily have $c=0$. When $U$ is complex, $c=n/2-\mu(E_{ss})$ where $n$ is the level of a Hodge representation.

\section{Example: Classification of Level 1 Hodge Representations}
In this section, we will classify the Hodge representations with Hodge numbers $h=(a,a)$ for some positive integer $a$. More specifically, we aim to classify the tuples $(\mf{g}_\IR,E, V_\IR)$, which according to Subsections 1.4-1.6 amounts to classifying tuples $(\mf{g}_\IC, E_{ss}, \mu,c)$ that corresponds to a level 1 Hodge representation as defined in Subsection 1.3. Namely $E$ acts on $V_\IC$ with eigenspace decomposition $V_\IC=V_{\frac{1}{2}}\oplus V_{-\frac{1}{2}}$, and $\dim V_\frac{1}{2}=\dim V_{-\frac{1}{2}}=a$. Partial results have been obtained by others, such as \cite{83paper}. 

\subsection{Reduction to Irreducible Representation}\label{Irrep}
Choose a grading element $E\in \mf{g}_\IC$. According to the discussion in Subsection 1.5, we can write $\mf{g}_\IC=\mf{g}_\IC^{ss}\oplus \mf{z}$, where $\mf{z}$ is a zero or one dimensional center, and $\mf{g}_\IC^{ss}=[\mf{g}_\IC,\mf{g}_\IC]$ is the semisimple part. Then we can write grading element $E=E_{ss}+E'$, where $E_{ss} \in \mf{g}_\IC^{ss}$ and $E'\in \mf{z}$. Assuming $V_\IR$ is irreducible, by Subsection 1.6, there are three possible cases: $U$ is a real, complex or quaternionic representation of $\mf{g}_\IR^{ss}$.
In case of level 1 Hodge representations, we want $$V_\IR\otimes \IC=V_{1/2}\oplus V_{-1/2}.$$ 

By discussion in Subsections 1.4-1.6, it suffices to locate tuples \\$(\mf{g}^{ss}_\IC,E_{ss},\mu,c)$, where $\mf{g}_\IC$ is a complex semisimple Lie algebra, $E_{ss}$ is a grading element of $\mf{g}^{ss}_\IC$, and $\mu$ is the highest weight on $\mf{g}_\IC^{ss}$'s complex irreducible representation $U$. We first claim that $\mf{g}^{ss}_\IC$ must be simple. To see this, let $\mf{g}_\IC^{ss}=\oplus_{i=1}^n \mf{g}_i$ where each $\mf{g}_i$ is simple. Then each irreducible representation of $\mf{g}_\IC^{ss}$ decomposes as $U=\otimes_{i=1}^n U_i$ where $U_i$ is an irreducible representation of $\mf{g}_i$. Moreover, suppose that $E$ is a grading element of $\mf{g}_\IC^{ss}$, then there exists a unique grading element $E_i$ from each $\mf{g}_i$ such that $E=\sum_{i=1}^n E_i$, and $\mu(E)+\mu^*(E)=\sum_{i=1}^n \mu_i(E_i)+\mu^*(E_i)$. Note that for each $i$, $\mu_i(E_i)+\mu^*_i(E_i)\geq 1$, and the tuples we search for must satisfy $\mu(E)+\mu^*(E)=1$. Hence, $\mf{g}_\IC^{ss}$ must be simple as claimed.

Now $U$ decomposes into $E_{ss}$'s eigenspaces $U_{\mu(E_{ss})}$ and $U_{\mu(E_{ss})-1}$, and correspondingly $U^*$ to decompose into $E_{ss}$'s eigenspaces $U^*_{1-\mu(E_{ss})}$ and $U^*_{-\mu(E_{ss})}$. Recall that $E'$ acts on $U$ via $c\mathds{1}$ and on $U^*$ via $-c\mathds{1}$ for some $c\in\IQ$. Then as $E=E_{ss}+E'$'s eigenspaces, $U$ and $U^*$ respectively decompose as:
$$U=U_{\mu(E_{ss})+c}\oplus U_{\mu(E_{ss})+c-1}$$
$$U^*=U_{1-\mu(E_{ss})-c}\oplus U^*_{-\mu(E_{ss})-c}$$
When $U$ is real, we must have $U=U_{1/2}\oplus U_{-1/2}$ as eigenspaces of $E_{ss}$, and $\dim U_{1/2}=\dim U_{-{1}/{2}}$. In this case, $c=0$. When $U$ is complex or quaternionic, with $c={1}/{2}-\mu(E_{ss})$, we get $U=U_{{1}/{2}}\oplus U_{-{1}/{2}}$ and $U^*=U^*_{{1}/{2}}\oplus U^*_{-{1}/{2}}$ as desired. Moreover, since $\dim U_{{1}/{2}}=\dim U^*_{-{1}/{2}}$ and $\dim U_{-{1}/{2}}=\dim U^*_{{1}/{2}}$, we have $\dim V_{1/2}=\dim V_{-1/2}$ as desired.

\begin{theorem}
The real irreducible Lie algebra Hodge representations with Hodge numbers ${{\bf{h}_\phi}}=(a,a)$ for some positive integer $a$ are given by the following tuples $(g_\IC^{ss},E_{ss},\mu,c)$:
\begin{enumerate}
    \item ($\mf{sl}(r+1,\IC),A^1,\omega_i,\frac{2i-r-1}{2(r+1)}$) for some $r\geq 1$, $1\leq i\leq r$, with \newline ${{{\bf{h}_\phi}}}=({\binom{r} {i-1}}+{\binom{r}{i}},{\binom{r} {i-1}}+{\binom{r}{i}})=({\binom{r}{i-1}},{\binom{r}{i}})+({\binom{r}{i}},{\binom{r}{i-1}})$ if $2i \neq r+1$, since in this case the representation is complex with respect to both $\mf{g}_{\IR}^{ss}$ and $\mf{g}_{\IR}$;
    \newline ${{{\bf{h}_\phi}}}=({\binom{r}{i-1}}+{\binom{r}{i}},{\binom{r}{i-1}}+{\binom{r}{i}})=({\binom{r}{i-1}},{\binom{r}{i}})+({\binom{r}{i}},{\binom{r}{i-1}})$ if $2i= r+1$ and $i$ is even, since in this case the representation is quaternionic with respect to both $\mf{g}_{\IR}^{ss}$ and $\mf{g}_{\IR}$;
    \\${{{\bf{h}_\phi}}}=(\frac{{\binom{r} {i-1}}+{\binom{r}{i}}}{2},\frac{{\binom{r}{i-1}}+{\binom{r}{i}}}{2})$ if $2i=r+1$ and $i$ is odd, since in this case the representation is real with respect to both $\mf{g}_{\IR}$ and $\mf{g}_{\IR}^{ss}$.
    \\ In all these cases, $\mf{g}_\IR^{ss}$ is $\mf{u}(1,r)$.
    \item ($\mf{sl}(r+1,\IC),A^r,\omega_i,\frac{r-2i+1}{2(r+1)}$) for some $r\geq 1$, $1\leq i\leq r$, with \newline ${{{\bf{h}_\phi}}}=(\binom{r}{i-1}+{\binom{r}{i}},{\binom{r}{i-1}}+{\binom{r}{i}})=(\binom{r}{i-1},\binom{r}{i})+({\binom{r}{i}},\binom{r}{i-1})$ if $2i \neq r+1$, since in this case the representation is complex with respect to both $\mf{g}_{\IR}$ and $\mf{g}_{\IR}^{ss}$;
    \newline ${{{\bf{h}_\phi}}}=(\binom{r} {i-1}+\binom{r}{i},\binom{r}{i-1}+\binom{r}{i})=(\binom{r}{i-1},\binom{r}{i})+(\binom{r}{i},\binom{r}{i-1})$ if $2i =r+1$ and $i$ is odd, since in this case the representation is quaternionic with respect to both $\mf{g}_{\IR}^{ss}$ and $\mf{g}_{\IR}$;
    \\${{{\bf{h}_\phi}}}=(\frac{\binom{r} {i-1}+\binom{r}{i}}{2},\frac{\binom{r} {i-1}+\binom{r}{i}}{2})$ if $2i=r+1$ and $i$ is even, since in this case the representation is real with respect to both $\mf{g}_{\IR}$ and $\mf{g}_{\IR}^{ss}$.
    \\ In all these cases, $\mf{g}_\IR^{ss}$ is $\mf{u}(r,1)$.
    
    \item ($\mf{sl}(r+1,\IC),A^i,\omega_1,\frac{2i-r-1}{2(r+1)}$) for some $r\geq 1$, $1\leq i\leq r$, with\newline ${{{\bf{h}_\phi}}}=(r+1,r+1)=(i,r+1-i)+(r+1-i,i)$ if $r\neq 1$, since in this case the representation is complex with respect to both $\mf{g}_{\IR}^{ss}$ and $\mf{g}_{\IR}$;
    \\ ${{{\bf{h}_\phi}}}=(1,1)$ if $r=1$, since in this case the representation is real with respect to both $\mf{g}_{\IR}^{ss}$ and $\mf{g}_{\IR}$.
    \\ In both cases, $\mf{g}_\IR^{ss}$ is $\mf{u}(i,r+1-i)$.
    
    \item ($\mf{sl}(r+1,\IC),A^i,\omega_r,\frac{r+1-2i}{2(r+1)}$) for some $r\geq 1$, $1\leq i\leq r$, with \newline ${{{\bf{h}_\phi}}}=(r+1,r+1)=(r+1-i,i)+(i,r+1-i)$ if $r\neq 1$, since in this case the representation is complex with respect to both $\mf{g}_{\IR}^{ss}$ and $\mf{g}_{\IR}$;
    \\ ${{\bf{h}_\phi}}=(1,1)$ if $r=1$, since in this case the representation is real with respect to both $\mf{g}_{\IR}^{ss}$ and $\mf{g}_{\IR}$.
     \\ In both cases, $\mf{g}_\IR^{ss}$ is $\mf{u}(i,r+1-i)$.
     
    \item ($\mf{so}(2r+1,\IC),A^1,\omega_r,0$) for some $r\geq 2$, with\newline ${{\bf{h}_\phi}}=(2^r,2^r)=(2^{r-1},2^{r-1})+(2^{r-1},2^{r-1})$ if $r\equiv 3,0\mod 4$, since in this case the representation is quaternionic with respect to both $\mf{g}_{\IR}^{ss}$ and $\mf{g}_{\IR}$;
    \\${{\bf{h}_\phi}}=(2^{r-1},2^{r-1})$ if $r\equiv 1,2\mod 4$, since in this case the representation is real with respect to both $\mf{g}_{\IR}^{ss}$ and $\mf{g}_{\IR}$.
     \\ In both cases, $\mf{g}_\IR^{ss}$ is $\mf{so}(2,2r-1)$.
     
    \item ($\mf{sp}(2r,\IC),A^{r},\omega_1,0$) with ${{\bf{h}_\phi}}=(r,r)$, since in this case the representation is real with respect to both $\mf{g}_{\IR}^{ss}$ and $\mf{g}_{\IR}$.
     In this case, $\mf{g}_\IR^{ss}$ is $\mf{sp}(r,\IR)$.
     
    \item ($\mf{so}(2r,\IC),A^{r-1}, \omega_1,0$) for some $r\geq 4$, with ${{\bf{h}_\phi}}=(2r,2r)=(r,r)+(r,r)$, since in this case the representation is quaternionic with respect to both $\mf{g}_{\IR}^{ss}$ and $\mf{g}_{\IR}$. In this case, $\mf{g}_\IR^{ss}$ is $\mf{so}^*(2r)$.
    
    \item ($\mf{so}(2r,\IC),A^{r}, \omega_1,0$ ) for some $r\geq 4$, with ${{\bf{h}_\phi}}=(2r,2r)=(r,r)+(r,r)$, since in this case the representation is quaternionic with respect to both $\mf{g}_{\IR}^{ss}$ and $\mf{g}_{\IR}$. In this case, $\mf{g}_\IR^{ss}$ is $\mf{so}^*(2r)$.
    
    \item ($\mf{so}(8,\IC),A^4, \omega_3,0$) with ${{\bf{h}_\phi}}=(8,8)=(4,4)+(4,4)$, since in this case the representation is quaternionic with respect to both $\mf{g}_{\IR}^{ss}$ and $\mf{g}_{\IR}$. In this case, $\mf{g}_\IR^{ss}$ is $\mf{so}^*(8)$.
    
    \item ($\mf{so}(8,\IC),A^1, \omega_3,0$) with ${{\bf{h}_\phi}}=(8,8)=(4,4)+(4,4)$, since in this case the representation is quaternionic with respect to both $\mf{g}_{\IR}^{ss}$ and $\mf{g}_{\IR}$. In this case, $\mf{g}_\IR^{ss}$ is $\mf{so}(2,6)$.
    
    \item ($\mf{so}(8,\IC),A^3, \omega_4,0$) with ${{\bf{h}_\phi}}=(8,8)=(4,4)+(4,4)$, since in this case the representation is quaternionic with respect to both $\mf{g}_{\IR}^{ss}$ and $\mf{g}_{\IR}$. In this case, $\mf{g}_\IR^{ss}$ is $\mf{so}^*(8)$.
    
    \item ($\mf{so}(8,\IC),A^1, \omega_4,0$) with ${{\bf{h}_\phi}}=(8,8)=(4,4)+(4,4)$, since in this case the representation is quaternionic with respect to both $\mf{g}_{\IR}^{ss}$ and $\mf{g}_{\IR}$. In this case, $\mf{g}_\IR^{ss}$ is $\mf{so}(2,6)$.
    
\end{enumerate}
\begin{proof}
Following from our discussion in Subsection $\ref{Irrep}$, to classify the irreducible Hodge representations with ${{\bf{h}_\phi}}=(a,a)$, it suffices to locate tuples $(\mf{g}_\IC^{ss},E_{ss},\mu)$, such that $U$ decomposes into $E_{ss}$'s eigenspaces $U_{\mu(E_{ss})}$ and $U_{\mu(E_{ss})-1}$, and correspondingly $U'$ decomposes into $E_{ss}$'s eigenspaces $U^*_{1-\mu(E_{ss})}$ and $U^*_{-\mu(E_{ss})}$, i.e. $\mu(E_{ss})+\mu^*(E_{ss})=1$. The results listed above can be easily verified by consulting Appendix $\ref{highest_lowest}$: Let's illustrate the calculations with $\mf{g}_\IC^{ss}=\mf{sp}(r)$ for some $r\geq 3$. Since we are considering Level 1 cases here, we need $(\mu+\mu^*)(E_{ss})=1$, so we need \begin{align*}
    \sum_{i=1}^{r-1}2\mu^i\{\alpha_1(E_{ss})+2\alpha_2(E_{ss})+\ldots+(i-1)\alpha_{i-1}
(E_{ss})+i(\alpha_i(E_{ss})+\alpha_{i+1}(E_{ss})+\\\ldots+\alpha_r(E_{ss})\}+\mu^r\{\alpha_1(E_{ss})+2\alpha_2(E_{ss})+\ldots+r\alpha_r(E_{ss})\}=1
\end{align*} where each $\mu^i$ is a nonnegative integer defined at the bottom of first page in Appendix \ref{highest_lowest}. Recall that we've also assumed each $\alpha_i(E_{ss})\in \{0,1\}$. Thus, we must have $\mu^i=0$ for all $1\leq i\leq r-1$ and $\mu^r=1$, since otherwise the fact that $E_{ss}\neq 0$ implies that the right hand side is greater than or equal to 2. Moreover, we must have $\alpha_1(E_{ss})=1$ and $\alpha_i(E_{ss})=0$ for all $1<i\leq r$, or again the right hand side is greater than or equal to 2. The judgments of whether the representations are real, quaternionic, or complex are also made by consulting Appendix A and employing the test given in Subsection 1.6.
\end{proof}
\end{theorem}
\begin{remark}
With fixed $r$, the representations given in  3\&4 are dual to each other and give rise to the same $V_\IC$; the representations given in 1\&2, 7\&8, and 9, 10, 11\&12 are equivalent up to an automorphism of the Dynkin diagram.
\end{remark}

\section{Example: Classification of Level 3 Hodge Representations}
\subsection{Representations of Simple Lie Algebras}
In this section, we assume that $\mf{g}_\IC^{ss}$ is simple. We treat the situation when $\mf{g}_\IC^{ss}$ is semisimple in Section 3.2.

\subsubsection{Reduction to Irreducible Representation}
For the Hodge representations that correspond to subdomains of the CY 3-folds, we want 
$$V_\IR\otimes \IC=V_{3/2}\oplus V_{1/2} \oplus V_{-1/2}\oplus V_{-3/2},$$
and 
$$\dim V_{3/2}=\dim V_{-3/2}=1.$$
And by virtue of the IPR relation on the Mumford-Tate domain first mentioned in Subsection 1.5, we may assume $\alpha_i(E_{ss})\in\{0,1\}$ for each simple root $\alpha_i$ \cite{Robles2014}.
Thus, to classify such Hodge representations, it suffices to find tuples $(\mf{g}_\IC, E_{ss}, \mu)$, where $\mu$ denote the highest weight on $U$ such that
\begin{enumerate}
    \item $\mu(E_{ss})+\mu^*(E_{ss})\in \{1,2,3\}$;
    \item If $\mu(E_{ss})+\mu^*(E_{ss})=3$, we need $\dim U_{\mu(E)}=\dim U_{-\mu^*(E)}=1$;
    If $\mu(E_{ss})+\mu^*(E_{ss})=1$ or $2$, we need $\dim U_{\mu(E)}=1$;
    \item $\alpha_i(E)\in\{0,1\}.$
\end{enumerate}

To see this, first consider the case when $\mu(E_{ss})+\mu^*(E_{ss})=1$. In this case, we only get a level 3 Hodge representation if $U$ is complex or quaternionic. If so, without loss of generality, we may assume that $V$ decomposes into $E$ eigenspaces $U_{3/2}$ and $U_{1/2}$, whereas $U^*$ decomposes into $E$ eigenspaces $U^*_{-1/2}$ and $U^*_{-3/2}$, with $c=3/2-\mu(E_{ss})$. Note that $V_{3/2}=U_{3/2}$, $V_{1/2}=U_{1/2}$, $V_{-1/2}=U^*_{-1/2}$ and $V_{-3/2}=U^*_{-3/2}$. Since $\dim U_{3/2}=\dim U^*_{-3/2}$ and $\dim U_{1/2}=\dim U^*_{-1/2}$, we must have $\dim V_{3/2}=\dim V_{-3/2}$ and $\dim V_{1/2}=\dim V_{-1/2}$. Then the condition 5 listed in Section 1 forces 
\begin{equation}\label{dim_1}
   \dim U_{3/2}=\dim U^*_{-3/2}=1
\end{equation}  
Write $E_{ss}=\sum_{i=1}^r \alpha_i(E_{ss}) A^i$ where each $\alpha_i(E_{ss})\in\{0,1\}$ and write $\mu|_{\mf{g}^{ss}_\IC}=\sum \mu^i \omega_i$ where each $\mu^i\in \IZ_{\geq 0}$. Then by a argument about parabolic subalgebra (see Appendix B), we get that the condition $\eqref{dim_1}$ is equivalent to 
\begin{equation}\label{dim_constraint}
    \mu^i\neq 0 \Rightarrow \alpha_i(E)=1.
\end{equation}
Therefore, when $\mu(E_{ss})+\mu^*(E_{ss})=1$, to classify the desired tuples, it suffices to classify complex or quaternionic irreducible ${g_\IC^{ss}}$ representation $U$ that satisfies $\eqref{dim_constraint}$ and admits $E_{ss}$ eigenspace decomposition $$U=U_{\mu(E_{ss})}\oplus U_{\mu(E_{ss})-1}.$$
\begin{proposition}
When $(\mu+\mu^*)(E_{ss})=1$, the irreducible Lie algebra Hodge representations with Hodge numbers ${\bf{h}_\phi}=(1,a,a,1)$ arise with the following tuples of $(\mf{g}_\IC^{ss},E_{ss},\mu,c)$:
\begin{enumerate}
    \item ($\mf{sl}(r+1,\IC),A^1,\omega_1,\frac{r+3}{2(r+1)}$) with ${{\bf{h}_\phi}}=(1,r,r,1)=(1,r,0,0)+(0,0,r,1)$ with $r\geq 1$, which is complex with respect to both $\mf{g}_{\IR}^{ss}$ and $\mf{g}_{\IR}$. When $r=1$, the representation is real with respect to $\mf{g}_{\IR}^{ss}$ but complex with respect to $\mf{g}_{\IR}$; otherwise, it is complex  with respect to both $\mf{g}_{\IR}^{ss}$ and $\mf{g}_{\IR}$. In this case, $\mf{g}_\IR^{ss}$ is $\mf{u}(1,r)$.
    
    \item ($\mf{sl}(r+1,\IC),A^r,\omega_r,\frac{r+3}{2(r+1)}$) with ${{\bf{h}_\phi}}=(1,r,r,1)$ with $r\geq 1$, which is complex with respect to $\mf{g}_{\IR}$. When $r=1$, the representation is real with respect to $\mf{g}_{\IR}^{ss}$; otherwise, it is complex with respect to $\mf{g}_{\IR}^{ss}$ In this case, $\mf{g}_\IR^{ss}$ is $\mf{u}(r,1)$.
\end{enumerate}
\begin{proof}
Again, by consulting Appendix A, one may enumerate all tuples that satisfy $(\mu+\mu^*)(E_{ss})=1$. We are not going to illustrate the computations again as they directly mirror the calculations done in the proof to Theorem 2.1. The reader only needs to bear in mind an extra point that now since we demand the Hodge numbers to be in form of $(1,a,a,1)$, whenever $\mu^i\neq 0$, we must have $\alpha_i(E)=1$. As elaborated in previous discussion, we must have that $U$ is complex or quaternionic, and that if $\mu$ contains $\omega_i$, then $E_{ss}=\sum \alpha_j(E_{ss})A^j$ must satisfy that $\alpha_i(E_{ss})\neq 0$, to get $V_\IC$ with desired Hodge numbers. 
One subtle point worth noticing is that when $r=1$, by the test in Subsection 1.6, the representation $U$ is real with respect to $\mf{g}_{\IR}^{ss}$. But with a nontrivial center that acts on all of $U$ with eigenvalue $1$ and acts on all of $U^*$ with eigenvalue $-1$, $U$ is complex with respect to $\mf{g}_{\IR}$, and therefore meet our eigenspace dimension requirements. The judgments of whether the representations are real, quaternionic, or complex are made by consulting appendix A and employing the test given in Subsection 1.6.
\end{proof}
\end{proposition}
\begin{remark}
One might notice that for any fixed $r$, the tuples in 1 and 2 are equivalent up to an automorphism of the Dynkin diagram. 
\end{remark}

\indent Now suppose $\mu(E_{ss})+\mu^*(E_
{ss})=2$. In this case, we only get a level 3 Hodge representation if $U$ is complex or quaternionic. If so, $U$ decomposes into $E$'s eigenspaces $U_{3/2}$, $U_{1/2}$ and $U_{-1/2}$, whereas $U^*$ decomposes into $E$'s eigenspaces $U^*_{1/2}$, $U^*_{-1/2}$ and $U^*_{-3/2}$ with $c=3/2-\mu(E_{ss})$. Moreover, $\dim U_{3/2}=\dim U^*_{-3/2}$, $\dim U_{1/2}=\dim U^*_{-1/2}$, and $\dim U_{-1/2}=\dim U^*_{1/2}$. Hence, we must have $\dim V_{3/2}=\dim V_{-3/2}$ and $\dim V_{1/2}=\dim V_{-1/2}$. Recall that for the Hodge representations that correspond to subdomains of the CY 3-folds, we must have $\dim V_{3/2}=\dim V_{-3/2}=1$, which is equivalent to condition $\eqref{dim_constraint}$. Therefore, when $\mu(E_{ss})+\mu^*(E_{ss})=2$, to classify the desired tuples, it suffices to classify complex or quaternionic irreducible ${g_\IC^{ss}}$ representation $U$ that satisfies $\eqref{dim_constraint}$ and admits $E_{ss}$ eigenspace decomposition $$U=U_{\mu(E_{ss})}\oplus U_{\mu(E_{ss})-1}\oplus U_{\mu(E_{ss})-2}.$$ 
\begin{proposition}
When $(\mu+\mu^*)(E_{ss})=2$, the irreducible Lie algebra Hodge representations with Hodge numbers ${\bf{h}_\phi}=(1,a,a,1)$ arise with the following tuples of $(\mf{g}_\IC^{ss},E_{ss},\mu,c)$:
\begin{enumerate}
    \item ($\mf{sl}(r+1,\IC),A^1,2\omega_1,\frac{-r+3}{2(r+1)}$) where $r\geq 1$ with  ${{\bf{h}_\phi}}=(1,(r+1)(r+2)/2-1,(r+1)(r+2)/2-1,1=(1,r,\frac{r(r+1)}{2},0)+(0,\frac{r(r+1)}{2},r,1)$. When $r=1$, the representation is real with respect to $\mf{g}_{\IR}^{ss}$ but complex with respect to $\mf{g}_{\IR}$; When $r\neq 1$, the representation is complex with respect to both $\mf{g}_{\IR}^{ss}$ and $\mf{g}_{\IR}$.  In this case, $\mf{g}_\IR^{ss}$ is $\mf{u}(1,r)$.
    
    \item ($\mf{sl}(r+1,\IC),A^r,2\omega_r,\frac{-r+3}{2(r+1)}$) where $r\geq 1$ with ${{\bf{h}_\phi}}=(1,(r+1)(r+2)/2-1,(r+1)(r+2)/2-1,1)=(1,r,\frac{r(r+1)}{2},0)+(0,\frac{r(r+1)}{2},r,1)$. When $r=1$, the representation is real with respect to $\mf{g}_{\IR}^{ss}$ but complex with respect to $\mf{g}_{\IR}$; When $r\neq 1$, the representation is complex with respect to both $\mf{g}_{\IR}^{ss}$ and $\mf{g}_{\IR}$. In this case, $\mf{g}_\IR^{ss}$ is $\mf{u}(r,1)$.
    
    \item ($\mf{sl}(r+1,\IC),A^2,\omega_2,\frac{-r+7}{2(r+1)}$), with ${{\bf{h}_\phi}}=(1,\frac{r(r+1)}{2}-1,\frac{r(r+1)}{2}-1,1)=(1,2(r-1),\frac{(r-1)(r-2)}{2},0)+(0,\frac{(r-1)(r-2)}{2},2(r-1),1)$. When $r=3$, the representation is real with respect to $\mf{g}_{\IR}^{ss}$ but complex with respect to $\mf{g}_{\IR}$; When $r\neq 1$, the representation is complex with respect to both $\mf{g}_{\IR}^{ss}$ and $\mf{g}_{\IR}$. In this case, $\mf{g}_\IR^{ss}$ is $\mf{u}(2,r-1)$.
    
    \item ($\mf{sl}(r+1,\IC),A^{r-1},\omega_{r-1},\frac{-r+7}{2(r+1)}$), with ${{\bf{h}_\phi}}=(1,\frac{r(r+1)}{2}-1,\frac{r(r+1)}{2}-1,1)=(1,2(r-1),\frac{(r-1)(r-2)}{2},0)+(0,\frac{(r-1)(r-2)}{2},2(r-1),1)$. When $r=3$, the representation is real with respect to $\mf{g}_{\IR}^{ss}$ but complex with respect to $\mf{g}_{\IR}$; When $r\neq 1$, the representation is complex with respect to both $\mf{g}_{\IR}^{ss}$ and $\mf{g}_{\IR}$. In this case, $\mf{g}_\IR^{ss}$ is $\mf{u}(r-1,2)$.
    
    \item ($\mf{sl}(r+1,\IC),A^{1}+A^{i},\omega_{1},\frac{3}{2}-\frac{2r-i+1}{r+1}$) for some $1\leq i\leq r$ with ${{\bf{h}_\phi}}=(1,r,r,1)=(1,i-1,r+1-i,0)+(0,r+1-i,i-1,1)$. The representation is complex with respect to both $\mf{g}_{\IR}^{ss}$ and $\mf{g}_{\IR}$. 
    
    \item ($\mf{sl}(r+1,\IC),A^{i}+A^{r},\omega_{r},\frac{3}{2}-\frac{i+r}{r+1}$) with ${{\bf{h}_\phi}}=(1,r,r,1)=(1,r-i,i,0)+(0,i,r-i,1)$. The representation is complex with respect to both $\mf{g}_{\IR}^{ss}$ and $\mf{g}_{\IR}$.
    
    \item ($\mf{so}(2r+1,\IC),A^1,\omega_{1},\frac{1}{2}$) with ${{\bf{h}_\phi}}=(1,2r-1,2r-1,1)=(1,2r-2,1,0)+(0,1,2r-2,1)$. The representation is real with respect to $\mf{g}_{\IR}^{ss}$ but complex with respect to $\mf{g}_{\IR}$. In this case, $\mf{g}_\IR^{ss}$ is $\mf{so}(2,2r-1)$.
    
    \item ($\mf{sp}(2r,\IC),A^1,\omega_{1},\frac{1}{2}$) with ${{\bf{h}_\phi}}=(1,2r-1,2r-1,1)=(1,2r-2,1,0)+(0,1,2r-2,1)$. This representation is quaternionic with respect to $\mf{g}_{\IR}^{ss}$ but complex with respect to $\mf{g}_{\IR}$. In this case, $\mf{g}_\IR^{ss}$ is $\mf{sp}(1,r-1)$.
    
    \item($\mf{so}(2r,\IC),A^1,\omega_{1},\frac{1}{2}$) with ${{\bf{h}_\phi}}=(1,2r-1,2r-1,1)=(1,2r-2,1,0)+(0,1,2r-2,1)$. This representation is real with respect to $\mf{g}_{\IR}^{ss}$ but complex with respect to $\mf{g}_{\IR}$. In this case, $\mf{g}_\IR^{ss}$ is $\mf{so}(2,2r-2)$.
    
    \item ($\mf{so}(8,\IC),A^3,\omega_{3},\frac{1}{2}$) with ${{\bf{h}_\phi}}=(1,7,7,1)=(1,6,1,0)+(0,1,6,1)$. This representation is real with respect to $\mf{g}_{\IR}^{ss}$ but complex with respect to $\mf{g}_{\IR}$. In this case, $\mf{g}_\IR^{ss}$ is $\mf{so}^*(8)$.
    
    \item($\mf{so}(8,\IC),A^4,\omega_{4},\frac{1}{2}$) with ${{\bf{h}_\phi}}=(1,7,7,1)=(1,6,1,0)+(0,1,6,1)$. This representation is real with respect to $\mf{g}_{\IR}^{ss}$ but complex with respect to $\mf{g}_{\IR}$. In this case, $\mf{g}_\IR^{ss}$ is $\mf{so}^*(8)$.
    
    \item ($\mf{so}(10,\IC),A^4,\omega_{4},\frac{1}{4}$) with ${{\bf{h}_\phi}}=(1,15,15,1)=(1,10,5,0)+(0,5,10,1)$. The representation is complex with respect to both $\mf{g}_{\IR}^{ss}$ and $\mf{g}_{\IR}$. In this case, $\mf{g}_\IR^{ss}$ is $\mf{so}^*(10)$.
    
    \item ($\mf{so}(10,\IC),A^5,\omega_{5},\frac{1}{4}$) with ${{\bf{h}_\phi}}=(1,15,15,1)=(1,10,5,0)+(0,5,10,1)$. The representation is complex with respect to both $\mf{g}_{\IR}^{ss}$ and $\mf{g}_{\IR}$. In this case, $\mf{g}_\IR^{ss}$ is $\mf{so}^*(2r)$.
    
    \item ($\mf{e}_6,A^1,\omega_{1},\frac{1}{6}$) with ${{\bf{h}_\phi}}=(1,26,26,1)=(1,16,10,0)+(0,10,16,1)$. The representation is complex with respect to both $\mf{g}_{\IR}^{ss}$ and $\mf{g}_{\IR}$. 
    
    \item ($\mf{e}_6,A^6,\omega_{6},\frac{1}{6}$) with ${{\bf{h}_\phi}}=(1,26,26,1)=(1,16,10,0)+(0,10,16,1)$. The representation is complex with respect to both $\mf{g}_{\IR}^{ss}$ and $\mf{g}_{\IR}$.
    
\end{enumerate}
\begin{proof}
By consulting Appendix A, one can enumerate all tuples that satisfy $(\mu+\mu^*)(E_{ss})=2$. As elaborated in the discussion preceding this proposition, to get $V_\IC$ with the desired Hodge numbers, we must have that $U$ is complex or quaternionic, and that if $\mu$ contains $\omega_i$ as a summand, then $E=\sum \alpha_j(E)A^j$ must satisfy $\alpha_i(E_{ss})\neq 0$.
For 1 and 2, when $r=1$ the representation is real with respect to $\mf{g}_{\IR}^{ss}$, but complex with respect to $\mf{g}_{\IR}$. For 3 and 4, when $r=4$ the representation is also real with respect to $\mf{g}_{\IR}^{ss}$, but complex with respect to $\mf{g}_{\IR}$. Similarly, ($\mf{so}(2r,\IC),A^1,\omega_{1},\frac{1}{2}$) and ($\mf{so}(2r+1,\IC),A^1,\omega_{1},\frac{1}{2}$) are real with respect to $\mf{g}_{\IR}^{ss}$ but complex with respect to $\mf{g}_{\IR}$, and so are ($\mf{so}(8,\IC),A^3,\omega_{3},\frac{1}{2}$) and ($\mf{so}(8,\IC),A^4,\omega_{4},\frac{1}{2}$). The reader might expect to see ($\mf{so}(5,\IC),A^2,\omega_{2},\frac{1}{2}$) with ${{\bf{h}_\phi}}=(1,3,3,1)$. This representation is quaternionic with respect to $\mf{g}_{\IR}^{ss}$ and complex with respect to $\mf{g}_{\IR}$, but it is equivalent to Item 7 with $r=2$, so we exclude it from the list.
\end{proof}
\end{proposition}
\begin{remark}
With fixed $r$, the representations given in 1\&2, 3\&4, 5\&6, 10\&11, 12\&13 and 14\&15 are equivalent up to an automorphism of the Dynkin diagram.
\end{remark}

\indent Finally, assume $\mu(E_{ss})+\mu^*(E_{ss})=3$. In this case, $U$ decomposes into $E_{ss}$ eigenspaces $U_{\mu(E_{ss})}$, $U_{\mu(E_{ss})-1}$, $U_{\mu(E_{ss})-2}$, and $U_{\mu(E_{ss})-3}$. Since $V_\IC$ decomposes into 4 eigenspaces of $E_{ss}$ with eigenvalues $\frac{3}{2},\frac{1}{2},-\frac{1}{2}$ and $-\frac{3}{2}$, we must have $\mu(E_{ss})=\frac{3}{2}$. Then the requirement that $\dim V_{3/2}=\dim V_{-3/2}=1$ implies that $U$ needs to be real and $\dim U_{3/2}=\dim U_{-3/2}=1$. Moreover, we need $c=0$ and equivalently $\mf{g}_\IC^{ss}=\mf{g}_\IC$. Since the dimension of the lowest weight space is necessarily equal to the dimension of the highest weight space, it suffices to require $\dim U_{3/2}=1$ and again this is equivalent to condition $\ref{dim_constraint}$. Therefore, when $\mu(E)+\mu^*(E)=3$, to classify the desired tuples, it suffices to classify real irreducible ${g_\IC^{ss}}$ representation $U$ that satisfies $\ref{dim_constraint}$ and admits $E_{ss}$ eigenspace decomposition $$U=U_{3/2}\oplus U_{1/2}\oplus U_{-1/2}\oplus U_{-3/2}.$$ 

\begin{proposition}
When $\mu+\mu^*(E)=3$, the irreducible Lie algebra Hodge representations with Hodge numbers ${\bf{h}_\phi}=(1,a,a,1)$ arise with the following tuples of $(\mf{g}_\IC^{ss},E_{ss},\mu,c)$:
\begin{enumerate}
    \item ($\mf{sl}(2,\IC),A^1,3\omega_1,0$) with ${{\bf{h}_\phi}}=(1,1,1,1)$. The representation is real with respect to both $\mf{g}_{\IR}^{ss}$ and $\mf{g}_{\IR}$. In this case, $\mf{g}_\IR^{ss}$ is $\mf{su}(1,r)$.
    
    \item ($\mf{sl}(6,\IC),A^3,\omega_3,0$) with ${{\bf{h}_\phi}}=(1,9,9,1)$.
    The representation is real with respect to both $\mf{g}_{\IR}^{ss}$ and $\mf{g}_{\IR}$. In this case, $\mf{g}_\IR^{ss}$ is $\mf{su}(3,3)$.
    
    \item ($\mf{so}(5,\IC),A^1+A^2,\omega_{2},0$) with ${{\bf{h}_\phi}}=(1,1,1,1)$.
    The representation is real with respect to both $\mf{g}_{\IR}^{ss}$ and $\mf{g}_{\IR}$. 
    
    \item ($\mf{sp}(6,\IC),A^3,\omega_{3},0$) with ${{\bf{h}_\phi}}=(1,6,6,1)$.
    The representation is real with respect to both $\mf{g}_{\IR}^{ss}$ and $\mf{g}_{\IR}$. In this case, $\mf{g}_\IR^{ss}$ is $\mf{sp}(6,\IR)$.
    
    \item ($\mf{sp}(2r,\IC),A^1+A^r,\omega_{1},0$) with ${{\bf{h}_\phi}}=(1,r-1,r-1,1)$.
    The representation is real with respect to both $\mf{g}_{\IR}^{ss}$ and $\mf{g}_{\IR}$. 
    
    \item ($\mf{so}(12,\IC),A^{5},\omega_{5},0$) with ${{\bf{h}_\phi}}=(1,15,15,1)$. The representation is real with respect to both $\mf{g}_{\IR}^{ss}$ and $\mf{g}_{\IR}$. In this case, $\mf{g}_\IR^{ss}$ is $\mf{so}^*(12)$.
    
    \item ($\mf{so}(12,\IC),A^{6},\omega_{6},0$) with ${{\bf{h}_\phi}}=(1,15,15,1)$. The representation is real with respect to both $\mf{g}_{\IR}^{ss}$ and $\mf{g}_{\IR}$. In this case, $\mf{g}_\IR^{ss}$ is $\mf{so}^*(12)$.
    
    \item ($\mf{e}_7,A^7,\omega_{7},0$) with ${{\bf{h}_\phi}}=(1,27,27,1)$. The representation is real with respect to both $\mf{g}_{\IR}^{ss}$ and $\mf{g}_{\IR}$. In this case, $\mf{g}_\IR^{ss}$ is $\mf{e}(6)\oplus\IR$.
\end{enumerate}
\begin{proof}
By consulting Appendix A, one can enumerate all tuples that satisfy $(\mu+\mu^*)(E)=3$. As elaborated in discussion preceding this proposition, to get $V_\IC$ with desired Hodge numbers, we must have that $U$ is real, and that if $\mu$ contains $\omega_i$ as a summand, then $E$ must contain $A^{i}$ as a summand.
Thus $(\mf{sl}(r+1,\IC),A^1+A^i+A^j,\omega_1)$, $(\mf{sl}(r+1,\IC),A^i+A^j+A^r,\omega_r)$,$(\mf{sl}(r+1,\IC),A^1+2A^i,\omega_1)$, $(\mf{sl}(r+1,\IC),A^r+2A^i,\omega_r)$,$(\mf{sl}(r+1,\IC),A^1+A^2,\omega_2)$,$(\mf{sl}(r+1,\IC),A^r+A^2,\omega_2)$,$(\mf{sl}(r+1,\IC),A^1+A^{r-1},\omega_{r-1})$, $(\mf{sl}(r+1,\IC),A^{r-1}+A^r,\omega_{r-1})$, ($\mf{so}(7,\IC),A^3,\omega_{3}$),($\mf{so}(10,\IC),A^1+A^5,\omega_{5}$),($\mf{so}(10,\IC),A^1+A^4,\omega_{4}$),($\mf{so}(10,\IC),A^4+A^5,\omega_{5}$),($\mf{so}(14,\IC),A^6,\omega_{6}$),($\mf{so}(14,\IC),A^7,\omega_{7}$), ($\mf{so}(8,\IC),A^1+A^{3},\omega_{3}$),\\ ($\mf{so}(8,\IC),A^3+A^{4},\omega_{3}$), ($\mf{so}(8,\IC),A^1+A^{4},\omega_{4}$), ($\mf{so}(8,\IC),A^3+A^{4},\omega_{4}$),\\ $(\mf{so}(2r,\IC), A^1+A^{r-1}, \omega_1)$ and $(\mf{so}(2r,\IC), A^1+A^{r}, \omega_1)$ are not on the list because they are either quaternionic or complex with respect to $\mf{g}_{\IR}^{ss}=\mf{g}_\IR$.
\end{proof}
\end{proposition}
\begin{remark}
With fixed $r$, the representations given in  6\&7 are equivalent up to an automorphism of the Dynkin diagram.
\end{remark}

\subsection{Representations of Semisimple Lie Algebras}
Now assume $\mf{g_\IC^{ss}}$ is semisimple. In this section, we will use the same notations as before. We claim that there are three possible cases in which we could get desired irreducible representation $V_\IC$ and grading element $E$ such that as eigenspaces of $E$, $V_\IC$ decomposes into $$V_\IC=V_{3/2}\oplus V_{1/2}\oplus V_{-1/2}\oplus V_{-3/2}$$ and the dimensions of $V_{3/2},V_{1/2},V_{-1/2},V_{-3/2}$ are respectively $1,a,a,1$ for some positive integer $a$. To see this, let $\mf{g_\IC^{ss}}=\oplus_{i=1}^n \mf{g}_i$ where each $\mf{g}_i$ is simple. Then each irreducible representation of $\mf{g_\IC^{ss}}$ decomposes as $U=\otimes_{i=1}^n U_i$ where $U_i$ is an irreducible representation of $\mf{g}_i$. Moreover, suppose that $E$ is a grading element of $\mf{g_\IC^{ss}}$, then there exists unique grading element $E_i$ from each $\mf{g}_i$ such that $E=\sum_{i=1}^n E_i$, and $\mu(E)+\mu^*(E)=\sum_{i=1}^n \mu_i(E_i)+\mu^*(E_i)$. Recall that for each $i$, $\mu_i(E_i)+\mu^*_i(E_i)\geq 1$, and the tuples we search for must satisfy $1\leq \mu(E)+\mu^*(E)\leq 3$. Thus, we conclude there are only three possible cases when $\mf{g_\IC^{ss}}$ is a direct sum of simple Lie algebras:
\begin{enumerate}
    \item $\mf{g_\IC^{ss}}=\mf{g}_1\oplus\mf{g}_2$ with irreducible representation $U=U_1\otimes U_2$ and $\mu_1(E_1)+\mu_1^*(E_1)=1$ and $\mu_2(E_2)+\mu_2^*(E_2)=1$;
    \item $\mf{g_\IC^{ss}}=\mf{g}_1\oplus\mf{g}_2$ with irreducible representation $U=U_1\otimes U_2$ and $\mu_1(E_1)+\mu_1^*(E_1)=1$ and $\mu_2(E_2)+\mu_2^*(E_2)=2$;
    \item $\mf{g_\IC^{ss}}=\mf{g}_1\oplus\mf{g}_2\oplus\mf{g}_3$ with irreducible representation $U=U_1\otimes U_2\otimes U_3$ and $\mu_1(E_1)+\mu_1^*(E_1)=\mu_2(E_2)+\mu_2^*(E_2)=\mu_3(E_3)+\mu_3^*(E_3)=1$.
\end{enumerate}

\subsubsection{Case 1: \texorpdfstring{$\mf{g_\IC^{ss}}$} is the direct sum of two simple Lie algebras, and \texorpdfstring{$U$} decomposes into three eigenspaces}
Suppose $\mf{g_\IC^{ss}}=\mf{g}_1\oplus \mf{g}_2$ where $\mf{g}_1$ and $\mf{g}_2$ are both simple. Then every grading element $E\in \mf{g_\IC}$ decomposes into  $$E=E_{ss}^1+E_{ss}^2+E',$$ 
where $E_{ss}^1\in \mf{g}_1$, $E_{ss}^2\in \mf{g}_2$ and $E'\in \mf{z}$. Even every irreducible representation $U$ of $\mf{g}_{\IC}^{ss}$ is isomorphic to the tensor product of some $U_1$ and $U_2$, where $U_1$ is an irreducible representation of $\mf{g}_1$ and $U_2$ is an irreducible representation of $\mf{g}_2$. In case 1, we assume that $U_1=U_a^1\oplus U_{a-1}^1$ as eigenspaces of $E_1$ and $U_2=U_b^2\oplus U_{b-1}^2$ as eigenspaces of $E_2$. Then $U_1\otimes U_2$ admits $E_{ss}^1+E_{ss}^2$ eigenspace decomposition
$$U=U_{a+b}\oplus U_{a+b-1} \oplus U_{a+b-2}.$$
With $c=3/2-(a+b)$ and $E'$ acting on $U$ as $c\mathds{1}$, we get that $U_1\otimes U_2$ admits $E$ eigenspace decomposition
$$U=U_{3/2}\oplus U_{1/2}\oplus U_{-1/2}.$$
Moreover, we have $\dim U_{3/2}=1$ since $\dim U_a=\dim U_b=1$. Thus, we get $V_\IC=U\oplus U^*$ a desired representation of $\mf{g}_\IC$ if $U$ is complex or quaternionic with respect to $\mf{g}_\IR$. Since the tensor product of two quaternionic representations is a real representation with respect to $\mf{g}_\IR$ by the test given in Subsection 1.6. 
There are only two possible cases:
\begin{enumerate}
    \item If one of $U_1$ and $U_2$ is complex with respect to $\mf{g}_\IR$, then $U_1\otimes U_2$ is complex with respect to $\mf{g}_\IR$; 
    \item If one of $U_1$ and $U_2$ is real and the other is quaternionic with respect to $\mf{g}_\IR$, then $U_1\otimes U_2$ is quaternionic with respect to $\mf{g}_\IR$.
\end{enumerate}
Examining the tuples listed in Proposition 3.3, one should see that the second case is impossible, because there is no tuple $(g_{ss},E_{ss},\mu)$ where $g_{ss}$ is simple such that the representation $U$ is quaternionic. 
We summarize all desired $(\mf{g_{\IC}^{ss}},E_{ss},\mu,c)$ tuples in the following proposition:
\begin{proposition}
In case 1, the irreducible Lie algebra Hodge representations with Hodge numbers ${\bf{h_\phi}}=(1,a,a,1)$ arise with the following tuples of $(\mf{g}_{1,\IC}^{ss},\mf{g}_{2,\IC}^{ss},E^1_{ss},E^2_{ss},\mu^1,\mu^2,c)$, where all $r_1\geq 1,r_2>1$:
\begin{enumerate}
    \item $(\mf{sl}(r_1+1,\IC), \mf{sl}(r_2+1,\IC),A^1,A^1,\omega_1,\omega_1,\frac{3}{2}-\frac{r_1}{r_1+1}-\frac{r_2}{r_2+1})$ with \\$h_\phi=(1,r_1+r_2+r_1r_2,r_1+r_2+r_1r_2,1)=(1,r_1+r_2,r_1r_2,0)+(0,r_1r_2,r_1+r_2,1)$. The representation is complex with respect to both $\mf{g}_\IC^{ss}$ and $\mf{g}_\IC$. In this case, $\mf{g}_\IR^{ss}$ is $\mf{u}(1,r_1)\oplus\mf{u}(1,r_2)$.
    
    \item $(\mf{sl}(r_1+1,\IC), \mf{sl}(r_2+1,\IC),A^1,A^r,\omega_1,\omega_r,\frac{3}{2}-\frac{r_1}{r_1+1}-\frac{r_2}{r_2+1})$ with \\$h_\phi=(1,r_1+r_2+r_1r_2,r_1+r_2+r_1r_2,1)=(1,r_1+r_2, r_1r_2,0)+(0,r_1r_2,r_1+r_2,1)$. The representation is complex with respect to both $\mf{g}_\IC^{ss}$ and $\mf{g}_\IC$. In this case, $\mf{g}_\IR^{ss}$ is $\mf{u}(1,r_1)\oplus\mf{u}(r_2,1)$.
    
    \item $(\mf{sl}(r_1+1,\IC), \mf{sl}(r_2+1,\IC),A^r,A^r,\omega_r,\omega_r,\frac{3}{2}-\frac{r_1}{r_1+1}-\frac{r_2}{r_2+1})$ with \\$h_\phi=(1,r_1+r_2+r_1r_2,r_1+r_2+r_1r_2,1)=(1,r_1+r_2, r_1r_2,0)+(0,r_1r_2,r_1+r_2,1)$. The representation is complex with respect to both $\mf{g}_\IC^{ss}$ and $\mf{g}_\IC$. In this case, $\mf{g}_\IR^{ss}$ is $\mf{u}(r_1,1)\oplus\mf{u}(r_2,1)$.
\end{enumerate}
\begin{proof}
The proof is similar to previous ones.
\end{proof}
\end{proposition}
\begin{remark}
 One might note that with fixed $r_1$ and $r_2$, 1, 2 and 3 are all equivalent to each other up to an automorphism of the Dynkin diagram.
\end{remark}

\subsubsection{Case 2: \texorpdfstring{$\mf{g}$} is the direct sum of 2 simple Lie algebras, and \texorpdfstring{$U$} decomposes into four eigenspaces}
Suppose $\mf{g}=\mf{g}_1\oplus \mf{g}_2$, where $\mf{g}_1$ and $\mf{g}_2$ are both simple Lie algebras. Similar to the previous case, grading element $E\in \mf{g}_{\IC}$ decomposes as $E=E_{ss}^1+E_{ss}^2+E'$, and irreducible representation $U$ of $\mf{g}_\IC$ must be isomorphic to some $U_1\otimes U_2$, where $U_1$ and $U_2$ are respectively irreducible representations of $\mf{g}_1$ and $\mf{g}_2$. In this case, we assume that $U_1=U_a^1\oplus U_{a-1}^1\oplus U_{a-2}^1$ and $U_2=U_b^2\oplus U_{b-1}^2$ respectively as eigenspaces of $E_{ss}^1$ and $E_{ss}^2$. Then as eigenspaces of $E_{ss}=E_{ss}^1+E_{ss}^2$, 
$$U=U_{a+b}\oplus U_{a+b-1}\oplus U_{a+b-2}\oplus U_{a+b-3}.$$
With $c=3/2-(a+b)$, we get the decomposition of $U$ as eigenspaces of $E$:
$$U=U_{3/2}\oplus U_{1/2}\oplus U_{-1/2}\oplus U_{-3/2}.$$
Hence, a desired tuple arises if and only if $U$ is real, if and only if either $U_1$ and $U_2$ are both real or both quaternionic. By examining the tuples examined in Proposition 3.3 and 3.5, one should see that the only the case in which both $U_1$ and $U_2$ are real is feasible.
We summarize all desired $(\mf{g}_\IC,E_{ss},\mu,c)$ tuples in the following proposition:
\begin{proposition}
In case 2, the irreducible Lie algebra Hodge representations with Hodge numbers ${\bf{h_\phi}}=(1,a,a,1)$ arise with the following tuples of $({\mf{g}}^{ss}_{1,\IC},{\mf{g}}^{ss}_{2,\IC},E^1_{ss},E^2_{ss},\mu^1,\mu^2,c)$:
\begin{enumerate}
    \item $(\mf{sl}(2,\IC), \mf{sl}(2,\IC),A^1,A^1,\omega_1,2\omega_1,0)$ with \\${\bf{h_\phi}}=(1,2,2,1)$. The representation is real. 
    \\In this case, $\mf{g}_\IR^{ss}$ is $\mf{su}(1,1)\oplus\mf{su}(1,1)$.
    
    \item $(\mf{sl}(2,\IC), \mf{sl}(4,\IC),A^1,A^2,\omega_1,\omega_2,0)$ with \\${\bf{h_\phi}}=(1,5,5,1)$. The representation is real. 
    \\In this case, $\mf{g}_\IR^{ss}$ is $\mf{su}(1,1)\oplus\mf{su}(2,2)$.
    
    \item $(\mf{sl}(2,\IC), \mf{so}(2r,\IC),A^1,A^1,\omega_1,\omega_1,0)$ with \\${\bf{h_\phi}}=(1,2r-1,2r-1,1)$. The representation is real.
    \\In this case, $\mf{g}_\IR^{ss}$ is $\mf{su}(1,1)\oplus\mf{so}(2,2r-2)$.
    
    \item $(\mf{sl}(2,\IC), \mf{so}(8,\IC),A^1,A^3,\omega_1,2\omega_3,0)$ with \\${\bf{h_\phi}}=(1,34,34,1)$. The representation is real.
    \\In this case, $\mf{g}_\IR^{ss}$ is $\mf{su}(1,1)\oplus\mf{so}^*(8)$.
    
    \item $(\mf{sl}(2,\IC), \mf{so}(8,\IC),A^1,A^4,\omega_1,\omega_4,0)$ with \\${\bf{h_\phi}}=(1,34,34,1)$. The representation is real.
    \\In this case, $\mf{g}_\IR^{ss}$ is $\mf{su}(1,1)\oplus\mf{so}^*(8)$.
    
     \item $(\mf{sl}(2,\IC), \mf{so}(2r+1,\IC),A^1,A^1,\omega_1,\omega_1,0)$ with \\${\bf{h_\phi}}=(1,2r,2r,1)$. The representation is real.
    \\In this case, $\mf{g}_\IR^{ss}$ is $\mf{su}(1,1)\oplus\mf{so}(2,2r-1)$.
\end{enumerate}
\begin{proof}
Again, the proof is similar to the ones before. 
\end{proof}
\end{proposition}
\begin{remark}
The representations given in 4\&5 are equivalent up to an automorphism of the Dynkin diagram.
\end{remark}

\subsubsection{Case 3: \texorpdfstring{$\mf{g}$} is the direct sum of 3 simple Lie algebras, and \texorpdfstring{$U$} decomposes into four eigenspaces}
Suppose $\mf{g}=\mf{g}_1\oplus \mf{g}_2\oplus \mf{g}_3$, where $\mf{g_1}$, $\mf{g_2}$ and $\mf{g}_3$ are all simple Lie algebras. Similar to the previous cases, grading element $E\in \mf{g}_{\IC}$ decomposes as $E=E_{ss}^1+E_{ss}^2+E_{ss}^3+E'$, and irreducible representation $U$ of $\mf{g}_\IC$ must be isomorphic to some $U_1\otimes U_2\otimes U_3$, where $U_1,U_2$ and $U_3$ are respectively irreducible representations of $\mf{g}_1,\mf{g}_2$ and $\mf{g}_3$. In this case, we assume that $U_1=U_a^1\oplus U_{a-1}^1$,$U_2=U_b^2\oplus U_{b-1}^2$, and $U_3=U_c^2\oplus U_{c-1}^2$ respectively as eigenspaces of $E_{ss}^1,E_{ss}^2$ and $E_{ss}^3$. Then as eigenspaces of $E_{ss}=E_{ss}^1+E_{ss}^2+E_{ss}^3$, 
$$U=U_{a+b+c}\oplus U_{a+b+c-1}\oplus U_{a+b+c-2}\oplus U_{a+b+c-3}.$$

With $c=3/2-(a+b)$ and $E'$ acting on $U$ as $c\mathds{1}$, we get that $U_1\otimes U_2$ admits $E$ eigenspace decomposition
$$U=U_{3/2}\oplus U_{1/2}\oplus U_{-1/2}\oplus U_{-3/2}.$$
Thus, we get desired $V_\IC$ in this case if and only if $U$ is real. Then there are two possible cases:
\begin{enumerate}
    \item All of $U_1$, $U_2$ and $U_3$ are real;
    \item Exactly two of them are quaternionic and one is real.
\end{enumerate} 
However, one may observe that all tuples examined in Proposition 3.3 are either complex or real, so the second case is impossible. We summarize all desired $(\mf{g_\IC}^{ss},E_{ss},\mu,c)$ tuples in the following proposition:
\begin{proposition}
In case 3, the irreducible Lie algebra Hodge representations with Hodge numbers ${\bf{h_\phi}}=(1,a,a,1)$ arise with the following tuples of \\ $(\mf{g}_{1,\IC}^{ss},\mf{g}^{ss}_{2,\IC},\mf{g}^{ss}_{3,\IC},E^1_{ss},E^2_{ss},E^3_{ss},\mu^1,\mu^2,\mu^3,c)$:
\begin{enumerate}
    \item $(\mf{sl}(2,\IC), \mf{sl}(2,\IC),\mf{sl}(2,\IC),A^1,A^1,A^1,\omega_1,\omega_1,\omega_1,0)$ with \\${\bf h_\phi}=(1,3,3,1)$. The representation is real.
    \\In this case, $\mf{g}_\IR^{ss}$ is $\mf{su}(1,1)\oplus\mf{su}(1,1)\oplus\mf{su}(1,1)$.
\end{enumerate}
\begin{proof}
Again, the proof is similar to the ones before.
\end{proof}
\end{proposition}

\newpage

\begin{appendices}
\section{Weight Difference in Terms of Simple Roots}\label{highest_lowest}

Fix a complex simple Lie Algebra $\mf{g}$ and an irreducible representation $V$ of $\mf{g}$. Denote $\mu$ as the highest weight on $V$. Fix Cartan subalgebra $\mf{h}$ of $g$ and choose $E\in \mf{h}$. Suppose $V$ decomposes into weight spaces $\{v_i\}$ of $\mf{h}$ with corresponding weights $\{\lambda_i\}$. Recall that we define the action of $\mf{h}$ on the dual space $V^*$ via:
\begin{equation*}
    \begin{split}
        g\cdot v_i^*(v_j)&:=v_i^*(-g\cdot v_j)\\
        &=v_i^*(-\lambda_j(g) v_j))\\
        &=-\lambda_j(g) v_i^*(v_j)
    \end{split}
\end{equation*}

Thus, if we denote the highest weight on $V^*$ as $\mu^*$, then the lowest weight on $V$ is just $-\mu^*$. Moreover, let $w_0$ be the longest element in the Weyl group of the root system of $\mf{g}$. Then we know that $\mu^*=-w_0 (\mu)$. On the other hand, denote the simple roots of $\mf{g}$ as $\Sigma$. Then $w_0(\Sigma)=-\Sigma$, so $-w_0$ defines an isometry on $\Sigma$ and thus gives an automorphism of the Dynkin diagram. Hence, by looking at all possible automorphisms of the Dynkin diagrams and find the automorphism that corresponds to the longest word in the Weyl group, we can get the formula of $w_0$. The results are exhibited below:

\begin{longtable}[c]{ |p{5.5cm}|p{3.5cm}|p{4cm}| }

 \hline
 \endhead
 
 \hline
 \endfoot
 
 \hline
 \endlastfoot
 
  Complex Simple Lie Algebra & Root System Type & $-w_0$  \\
  \hline
   $\mf{sl}({r+1},\IC)$ & $A_r$ & $\omega_i\leftrightarrow\omega_{r+1-i}$  \\
   
   $\mf{so}(2r+1,\IC)$  & $B_r$ & id \\
   
   $\mf{sp}(2r,\IC)$ & $C_r$ & id   \\
   
   $\mf{so}(2r,\IC)$ with $r$ even & $D_r$ & id  \\
   $\mf{so}(2r,\IC)$ with $r$ odd & $D_r$ & $\omega_{r-1}\leftrightarrow \omega_r$  \\
   $\mf{e_6}$ & $E_6$ & $\omega_1\leftrightarrow \omega_6,\omega_3\leftrightarrow \omega_5$ \\
   $\mf{e_7}$ & $E_7$ & id \\
   $\mf{e_8}$ & $E_8$ & id \\
   $\mf{f_4}$ & $F_4$ & id \\
   $\mf{g_2}$ & $G_2$ & id \\
 \hline
\end{longtable}

Finally, by making use of the formulae that translate fundamental weights into simple roots, we may express the difference of the highest weight and the lowest weight on $V$, $\mu+\mu^*$ with simple roots. Write $\mu=\sum_{i=1}^r \mu^i\omega_i$, where $\omega_i$ denotes the $i$-th fundamental weight. The results are exhibited below:
\begin{enumerate}
    \item When $\mf{g}=\mf{sl}(r+1,\IC)$, 
    \begin{flalign*}
    \mu+\mu^*=\sum_{i=1}^{r}\mu^i\{\sum_{l=1}^i l\alpha_l+m\sum_{k=i+1}^{r-i}\alpha_k+\sum_{l=1}^i l\alpha_{r+1-l}\}
    \end{flalign*}
    where $m=\min\{i,r+1-i\}.$

    \item When $\mf{g}=\mf{so}(2r+1,\IC)$, 
    \begin{flalign*}
    \mu+\mu^* &=\sum_{i=1}^{r-1}2\mu^i\{\alpha_1+2\alpha_2+\ldots+(i-1)\alpha_{i-1}+i(\alpha_i+\alpha_{i+1}+\ldots+\alpha_{r})\}\\
    &+\mu^r\{\alpha_1+2\alpha_2+\ldots+r\alpha_r\}.
    \end{flalign*}
    
    \item When $\mf{g}=\mf{sp}(2r,\IC)$, 
    \begin{flalign*} 
    \mu+\mu^*&=\sum_{i=1}^{r}2\mu^i\{\alpha_1+2\alpha_2+\ldots+(i-1)\alpha_{i-1}\\
    &+i(\alpha_i+\alpha_{i+1}+\ldots+\alpha_{r-1}+\frac 1 2 \alpha_r)\}.
    \end{flalign*}
    
    \item When $\mf{g}=\mf{so}({2r},\IC)$ with $r$ odd,\begin{flalign*} 
    \mu+\mu^* &= \sum_{i=1}^{r-2} 2\mu^i\{ \alpha_1+2\alpha_2+\ldots+(i-1)\alpha_{i-1}\\
    &+i(\alpha_i+\alpha_{i+1}+\ldots+\alpha_{r-2})+\frac 1 2 i(\alpha_{r-1}+\alpha_r)\}\\
    &+\frac {\mu^{r-1}+\mu^r}{2}\{\alpha_1+2\alpha_2+\ldots+(r-2)\alpha_{r-2}+\frac r 2\alpha_{r-1}+\frac {r-2} 2 \alpha_r\}\\
    &+\frac{\mu^r+\mu^{r-1}}{2}\{\alpha_1+2\alpha_2+\ldots+(r-2)\alpha_{r-2}+\frac {r-2} 2 \alpha_{r-1}+\frac r 2 \alpha_r\}.
    \end{flalign*}
    
    \item When $\mf{g}=\mf{so}({2r},\IC)$ with $r$ even,
    \begin{flalign*} 
    \mu+\mu^* &= \sum_{i=1}^{r-2} 2\mu^i\{ \alpha_1+2\alpha_2+\ldots+(i-1)\alpha_{i-1}\\
    &+i(\alpha_i+\alpha_{i+1}+\ldots+\alpha_{r-2})+\frac 1 2 i(\alpha_{r-1}+\alpha_r)\}\\
    &+\mu^{r-1}\{\alpha_1+2\alpha_2+\ldots+(r-2)\alpha_{r-2}+\frac r 2\alpha_{r-1}+\frac {r-2} 2 \alpha_r\}\\
    &+\mu^r\{\alpha_1+2\alpha_2+\ldots+(r-2)\alpha_{r-2}+\frac {r-2} 2 \alpha_{r-1}+\frac r 2 \alpha_r\}.
    \end{flalign*}
    
    \item When $\mf{g}=\mf{e_6}$, 
    \begin{flalign*}
        \mu+\mu^*&= (\mu^1+\mu^6)(2\alpha_1+2\alpha_2+3\alpha_3+4 \alpha_4+3\alpha_5+2\alpha_6)\\
        &+(\mu^3+\mu^5)(3\alpha_1+4\alpha_2+6\alpha_3+8\alpha_4+6\alpha_5+3\alpha_6)\\
        &+2\mu^2(\alpha_1+2\alpha_2+2\alpha_3+3\alpha_4+2\alpha_5+\alpha_6)\\
        &+2\mu^4(2\alpha_1+3\alpha_2+4\alpha_3+6\alpha_4+4\alpha_5+2\alpha_6).
    \end{flalign*}
    
    \item When $\mf{g}=\mf{e_7}$, 
    \begin{flalign*}
        \mu+\mu^*&= 2\mu^1(2\alpha_1+2\alpha_2+3\alpha_3+4 \alpha_4+3\alpha_5+2\alpha_6+\alpha_7)\\
        &+\mu^2(4\alpha_1+7\alpha_2+8\alpha_3+12\alpha_4+9\alpha_5+6\alpha_6+3\alpha_7)\\
        &+2\mu^3(3\alpha_1+4\alpha_2+6\alpha_3+8\alpha_4+6\alpha_5+4\alpha_6+2\alpha_7)\\
        &+2\mu^4(4\alpha_1+6\alpha_2+8\alpha_3+12\alpha_4+9\alpha_5+6\alpha_6+3\alpha_7)\\
        &+\mu^5(6\alpha_1+9\alpha_2+12\alpha_3+18\alpha_4+15\alpha_5+10\alpha_6+5\alpha_7)\\
        &+2\mu^6(2\alpha_1+3\alpha_2+4\alpha_3+6\alpha_4+5\alpha_5+4\alpha_6+2\alpha_7)\\
        &+\mu^7(2\alpha_1+3\alpha_2+4\alpha_3+ 6\alpha_4+5\alpha_5+4\alpha_6+3\alpha_7).
    \end{flalign*}
    
    \item When $\mf{g}=\mf{e_8}$, 
    \begin{flalign*}
        \mu+\mu^*&= 2\mu^1(4\alpha_1+5\alpha_2+7\alpha_3+10 \alpha_4+8\alpha_5+6\alpha_6+4\alpha_7+2\alpha_8)\\
        &+2\mu^2(5\alpha_1+8\alpha_2+10\alpha_3+15\alpha_4+12\alpha_5+9\alpha_6+6\alpha_7+3\alpha_8)\\
        &+2\mu^3(7\alpha_1+10\alpha_2+14\alpha_3+20\alpha_4+16\alpha_5+12\alpha_6+8\alpha_7+4\alpha_8)\\
        &+2\mu^4(10\alpha_1+15\alpha_2+20\alpha_3+30\alpha_4+24\alpha_5+18\alpha_6+12\alpha_7+6\alpha_8)\\
        &+2\mu^5(8\alpha_1+12\alpha_2+16\alpha_3+24\alpha_4+20\alpha_5+15\alpha_6+10\alpha_7+5\alpha_8)\\
        &+2\mu^6(6\alpha_1+9\alpha_2+12\alpha_3+18\alpha_4+15\alpha_5+12\alpha_6+8\alpha_7+4\alpha_8)\\
        &+2\mu^7(4\alpha_1+6\alpha_2+8\alpha_3+ 12\alpha_4+10\alpha_5+8\alpha_6+6\alpha_7+3\alpha_8)\\
        &+2\mu^8(2\alpha_1+3\alpha_2+4\alpha_3+ 6\alpha_4+5\alpha_5+4\alpha_6+3\alpha_7+2\alpha_8).
    \end{flalign*}
    
    \item When $\mf{g}=\mf{f_4}$, 
    \begin{flalign*}
        \mu+\mu^*&= 2\mu^1(2\alpha_1+3\alpha_2+4\alpha_3+2 \alpha_4)
        +2\mu^2(3\alpha_1+6\alpha_2+8\alpha_3+4\alpha_4)\\
        &+2\mu^3(2\alpha_1+4\alpha_2+6\alpha_3+3\alpha_4)
        +2\mu^4(\alpha_1+2\alpha_2+3\alpha_3+2\alpha_4).
    \end{flalign*}
    
    \item When $\mf{g}=\mf{g_2}$, 
    \begin{flalign*}
        \mu+\mu^*&= 2\mu^1(2\alpha_1+\alpha_2)
        +2\mu^2(3\alpha_1+2\alpha_2).    \end{flalign*}
\end{enumerate}

\section{Parabolic Subalgebra Proof}
\label{parabolic}
In this section, we will define Borel subalgebra and parabolic subalgebra. Then we will prove the statement that given a  complex semisimple Lie algebra $\mf{g}_\IC$, its irreducible representation $U$ with highest weight $\mu$ and grading element $E\in\mf{g}_\IC$, $\dim U_\mu(E)=1$ if and only if $\alpha_i(E)>0$ for all $\mu^i\neq 0$.  

\begin{definition}
Given $\mf{g^{ss}}$ a complex semisimple Lie algebra, $\mf{h}\subset\mf{g}^{ss}$ a Cartan subalgebra, and $\Delta^{+}\subset \mf{h}^*$ a set of positive roots. The \textbf{Borel subalgebra} determined by $(\mf{h},\Delta^+)$ is 
$$\mf{b}=\mf{h}\oplus \bigoplus_{\alpha\in\Delta^+}\mf{g}_{\alpha}.$$
\end{definition}

\begin{definition}
A \textbf{parabolic subalgebra} is a subalgebra that contains a Borel subalgebra. 
\end{definition}

\begin{remark}
Given $\mf{g^{ss}}$ a complex semisimple Lie algebra of rank $r$, $\mf{h}\subset\mf{g}^{ss}$ a Cartan subalgebra, $\Delta^{+}\subset \mf{h}^*$ a set of positive roots, and $I\subset\{1,\ldots,r\}$, the {parabolic subalgebra} determined by $(\mf{h},\Delta^{
+},I)$ is a subalgebra that admits decomposition:
$$\mf{p}_I=\mf{h}\oplus \bigoplus_{\alpha\in\Delta_I}\mf{g}_{\alpha}\oplus\bigoplus_{\alpha\in\Delta_I^+}\mf{g}_\alpha,$$
where $\Delta_I=\{\sum_{i=1}^r \lambda^i\alpha_i | \lambda^i=0 \forall i\in I\}$ and $\Delta_I^+=\Delta^+\setminus\Delta_I.$ In fact, given Borel subalgebra $\mf{b}$ determined by $(\mf{h},\Delta^+)$, each parabolic subalgebra $\mf{p}\supset \mf{b}$ is determined by a unique maximal set $I$. 
\end{remark}
Now fix a complex semisimple Lie algebra $\mf{g}_\IC$. Let $U$ be an irreducible representation of $\mf{g}_\IC$ with highest weight $\mu=\sum_{i=1}^r \mu^i \omega_i$ and $0\neq v\in U_\mu$ a highest weight vector.

\begin{proposition}
$\mf{p}_\mu:=\{\xi\in\mf{g_\IC}|\xi\cdot v\in \IC v\}$ is a parabolic subalgebra.
\begin{proof}
Since all elements of the Cartan subalgebra act on $U$ via weights, $\mf{h}\subset \mf{p}_I$. Since $v$ is a highest weight vector, it is annihilated by all positive root vectors, so all positive root spaces are contained in $\mf{p}_\mu$.
 Thus, $\mf{p}_\mu$ contains a Borel subalgebra and is therefore a parabolic subalgebra.
\end{proof}
\end{proposition}

\begin{proposition}
$\mf{p}_\mu$ is the subalgebra determined by $(\mf{h},\Delta^+,I)$ where $\Delta^+$ is the set of all positive roots and $I=\{i|\mu^i\neq 0\}$. 
\begin{proof}
Let $I'$ be the maximal set that determines $\mf{p}_\mu$. It is straightforward that $I\subset I'$. We will show that $I'\subset I$. Now suppose that there exists some $i$ such that in the sum $\mu=\sum_{j=1}^r \mu^j\omega_j$, $\mu^i\neq 0$ and $\mf{p}_\mu$ contains some $0\neq Y\in \mf{g}_{-\alpha_i}$. We first claim that $Y\cdot v=0$. Since $Y\in\mf{p}_\mu$, $Y\cdot v=cv$ for some $c\in \IC$. Thus, for all $H\in \mf{h}$, $$H\cdot Y\cdot v=cH\cdot v=c\mu(H)v.$$
On the other hand, $$H\cdot Y\cdot v=([H,Y]+Y\cdot H)\cdot v=(-\alpha_i(H)+1)Y\cdot v=(-\alpha_i(H)+1)cv.$$ Hence, we must have $c=0$.
\\ Choose $X\in \mf{g}_{\alpha_i}\subset \mf{p}_\mu$. Since $v$ is a highest weight vector, $X\cdot v=0$. Define $H:=[X,Y]$. Then $\omega_i(H)> 0$ and thus $\mu(H)>0$. However,
$$\mu(H)v=H\cdot v=[X,Y]\cdot v= X\cdot Y\cdot v-Y\cdot X\cdot v=X\cdot 0-Y\cdot 0 =0;$$ contradiction! Thus, we must have I'=I.
\end{proof}
\end{proposition}

\begin{corollary}
Suppose $0\neq \xi\in\mf{g}_{-\alpha_i}$. Then $\xi\cdot v\neq 0$ if and only if $i\in I$.
\begin{proof}
If $i\in I$, then  $\mf{g}_{-\alpha_i}\cap \mf{g}_\mu=\{0\}$. Thus, $\xi\cdot v \notin \IC v$ and thus $\xi\cdot v\neq 0$. If $\xi\circ v=0$, then $\xi\in \mf{g_\mu}$. More specifically, $\xi \in \bigoplus_{\alpha\in \Delta_I}\mf{g}_\alpha$. Thus, $i\in I$.
\end{proof}
\end{corollary}

\begin{proposition}
$\dim U_{\mu(E)}=1$ if and only if $\alpha_i(E)>0$ for all $i\in I$.
\begin{proof}
Fix $i\in I$ and choose $0\neq \xi\in \mf{g}_{-\alpha_i}$. Again, since $i\in I$, we have that $\mf{g}_{-\alpha_i}\cap \mf{g}_\mu =\{0\}$, so $\xi\notin \mf{g}_\mu$. Therefore, $\xi\cdot v\notin \IC v$. From Corollary B.6, we also know that $\xi\cdot v\neq 0$. Now note that
$$E\cdot \xi\cdot v=([E,\xi]+\xi\cdot E)\cdot v=(-\alpha_i(E)+\mu(E))\xi\cdot v.$$
Thus, $\alpha_i(E)=0$ if and only if $0\neq \xi\cdot v\in U_{\mu(E)}$, if and only if $\dim U_{\mu(E)}>1$. Hence, $\dim U_{\mu(E)}=1$ if and only if $\alpha_i(E)>0$ for all $i\in I$ as claimed.
\end{proof}
\end{proposition}
\section{Center Dimension Proof}
\label{center}
Given a complex Lie algebra $\mf{g}$, let $G$ be a complex connected Lie group with Lie algebra $\mf{g}$. We will prove that the center $Z(\mf{g})$ of $\mf{g}$ has at most dimension one if $\mf{g}/Z(\mf{g})$ is simple, so the adjoint representation of $\mf{g}$ is irreducible. For sake of contradiction, suppose that $\dim Z(\mf{g})\geq 2$. Recall that by definition $$Z(\mf{g})=\{X\in \mf{g}|[X,Y]=0 \quad\forall Y\in \mf{g}\}.$$ Then we can find basis $\{X_1,X_2,\ldots,X_m\}$ of $Z(\mf{g})$, where $m\geq 2$ is the dimension of $Z(\mf{g})$. Define a linear map $\phi:\mf{g}\to\mf{g}$ that switches $X_1$ and $X_2$ and fixes all other vectors in $\mf{g}$. Define the $G$ representation $Ad$ on $\mf{g}$ via $$Ad(exp(X))\cdot Y=[X,Y].$$ Note that by the connectedness of $G$ the above definition is well-defined. Now to show that $Ad$ commutes with $\phi$, it suffices to observe that: 
\begin{align*}
    \phi \circ Ad(exp(X))\cdot X_1 &=\phi [X,X_1]=0=\phi[X,X_2]\\&=Ad(exp(X))\cdot X_2 = Ad(exp(X))\cdot \phi(X_1)
\end{align*}
By Schurs's lemma, on a $G$-irreducible representation, the only $G$-linear automorphism must be a multiple of identity, but $\phi$ is not; contradiction. Hence, the dimension of any complex Lie algebra's center is at most one.
\end{appendices}

\medskip
\printbibliography
\end{document}